\documentclass[a4paper,10pt]{amsart}

\newcommand{\R}{\mathbf{R}}
\newcommand{\C}{\mathbf{C}}

\newtheorem{theorem}{Theorem}
\newtheorem{proposition}{Proposition}
\newtheorem{corollary}{Corollary}

\title{Dilational Hilbert Scales \\
and Deconvolutional Sharpening}
\author{}

\begin{document}

\maketitle

\begin{center}
M. Hegland, Centre for Mathematics and its Applications, ANU, Canberra, ACT~0200, Australia  
and R. S. Anderssen, 
CSIRO Mathematical and Information Sciences, PO Box 664, Canberra, ACT~2601, Australia
\end{center}

\begin{abstract}
Operationally, index functions of variable Hilbert scales can be viewed as generators 
for families of spaces and norms. 
Using a one parameter family of index functions based on the dilations of a
given index function, 
a new class of scales (dilational Hilbert scales (DHS)) 
is derived which generates new interpolatory inequalities 
(dilational interpolatory inequalities (DII))
which have the ordinary Hilbert scales (OHS) 
interpolatory inequalities as special cases.
They therefore represent a one-parameter family generalization of OHS,
and are a precise and concise subset of VHS approriate for deriving
error estimates for deconvolution.
The role of the Hilbert scales in deriving error estimates for 
the approximate solution of inverse problems is discussed along 
with an application of DHS to deconvolution sharpening.
\end{abstract}

\section{Introduction}

In the analysis of the numerical performance of traditional
regularization methods (Engl \emph{et al.}~\cite{EngHN96}),
interpolatory inequalities between the norms $||u||_\alpha$, $\alpha\in R_+$, 
generated by an appropriate family of Hilbert spaces
$H_\alpha=R(T^{\alpha/2})$, play the central 
role.\footnote{Here $R(T^{\alpha/2})$ denotes the range
of the operator $T^{\alpha/2}$.}
In terms of the original concept of a Hilbert scale, 
as introduced by Krein and Petunin~\cite{KreP66}
and generated by a densely-defined, unbounded, self-adjoint and 
strictly positive operator $T$,
such inequalities take,
for a given value $\gamma$ of the linear scale between given values $\alpha$ and $\beta$,
the form
\begin{equation}
  \|u\|_{\theta\alpha + (1-\theta)\beta} \leq \|u\|_\alpha^\theta\, \|u\|_\beta^{1-\theta}, \quad
  \text{for all $\theta\in (0,1)$ and $0 \leq \alpha < \beta$}
\end{equation}
where $(\cdot,\cdot)$ and $\|\cdot\|$ denote the inner product and norm in the original
Hilbert space $H=H_0$ and where 
$$
  \|u\|_\alpha = \|T^{\alpha/2} u \|, \quad \alpha > 0
$$
and similar for $\|u\|_\beta$.

Through the appropriate choice of $T$ and the values for $a$ and $b$, 
the corresponding inequality (1) can be used to derive error estimates
for the regularized solution of improperly posed operator equations
that simultaneously take account of both the compact and smoothing nature of
the operator 
(Groetsch~\cite{Gro84}, Natterer~\cite{Nat84}, 
Schr{\"o}ter and Tautenhahn~\cite{SchT94},
Tautenhahn~\cite{Tau96}).

When utilizing such inequalities to derive error estimates for linear
improperly posed problems,
it was observed by various authors that realistic
error estimates could only be derived for a subset of 
linear improperly posed operator equations.
This led to the need to construct more general 
counterparts of the inequality (1).
Hegland~\cite{Heg92},~\cite{Heg95}, 
by exploiting the spectral decomposition of an appropriately chosen 
operator $T$, 
first introduced the concept of a \emph{variable Hilbert scale} 
(VHS) for a quite general index function $\phi$.
Then, by invoking appropriate regularity about the choice of the index
function,
Hegland established how more general counterparts of the interpolatory
inequality (1) could be constructed.

The utility of this basic concept of a VHS,
in deriving more representative error estimates and convergence rates
for the regularized solution of improperly posed operator equations,
has been subsequently exploited and/or modified by various authors
including Nair \emph{et al.}~\cite{NaiPT05}, 
Math{\'e} and Pereverez~\cite{MatP03},~\cite{MatP06},
Hofmann and Yamamoto~\cite{HofY05}, 
Math{\'e} and Tautenhahn~\cite{MatT06} and B\'egout and Soria~\cite{BegS06}.

The concept of a variable Hilbert scale not only overcomes the mentioned 
shortcoming of interpolatory inequalities of the form (1),
but also allows a greater variety of interpolatory inequalities to be
constructed.
It is this aspect that is pursued in this paper.
A recent example can be found in Nair \emph{et al.}~\cite{NaiPT05}. 
In our paper, the concept of a family of variable Hilbert scales and associated index
functions are introduced using the spectral theorem for positive
definite self-adjoint operators.
The special subclass of dilational Hilbert scales (DHS) 
is then defined,
and a new type of interpolatory inequality is derived and applied.
 
The motivation is the analysis of deconvolution by sharpening (Hegland
and Anderssen~\cite{HegA05}).
When the Gaussian model 
\[
G_\gamma(t)=\frac{1}{\gamma}G(t/\gamma), \quad
G(t)=\exp(- t^2/2), 
\]
defines the smoothing kernel $B_\gamma$ 
(e.g. the broadening function that models the measurement process in
spectroscopy as examined in Section 3.3) in the convolution equation
\begin{equation}
\label{eq:Bfeqg}
 g=B_\gamma*f,
\end{equation}
the convolutional relationship 
\[
G_\alpha*G_\beta = G_\gamma, \quad
\gamma=\frac{\alpha\beta}{\sqrt{\alpha^2+\beta^2}},
\]
can be used to factor 
\[
g=G_\gamma*f
\]
in the following manner
\[
g=G_\alpha*G_\beta*f=
G_\alpha*z, \quad 
z=G_\beta*f
\]
with $\alpha$ determining the nature of the sharpening to be performed
and $\beta$ the form of the source condition.
In this way, $\alpha$ and $\beta$ perform a trade-off, 
mentioned earlier,
between the sharpening and the achievable rates of convergence.
Though such a strategy does not in general hold for non-Gaussian models
of the smoothing kernels, 
it can be adapted to dilational models
\[
B_\gamma(t) = \frac{1}{\gamma} B(t/\gamma), \quad B(t)\sim 
\text{some~representative~peak~function,}
\]
by performing the factorization in the following manner
\begin{equation}
g=B_{\gamma,\beta}*B_\beta*f= B_{\gamma,\beta}*z, 
\quad B_\gamma=B_{\gamma,\beta}*B_\beta, 
\quad z=B_\beta*f,
\end{equation}
where the form taken by $B_{\gamma,\beta}(t)$ depends on $B_\gamma$,
chosen to reflect the convolutional smoothing that the measurement
system has performed on broadening $f$ to give $g$, and on $\beta$,
chosen to ensure that the source condition encapsulated in $z=B_\beta*f$
yields acceptable convergence rates. A discussion about how $y=B_\beta*f$
provides an appropriate source condition is given in Section~\ref{sec:3.3}.

Now, however, it is necessary to determine and work with 
$B_{\gamma,\beta}(t)$
as the sharpening to be utilized when performing deconvolution.
The value of $\gamma$ reflects the difference in instruments with small
$\gamma$ corresponding to expensive accurate measurement instrumentation.
The scaling $1/\gamma$, introduced into the definition of $B_\gamma$, is
an essential regularity condition. As $\beta$ tends to zero, $B_\beta(t)$
must have a behaviour like that of the identity operator so that, under
appropriate circumstances, $z$ recovers $f$. For the scaling chosen here,
it guarantees that $z=f$ when $f$ is a constant.

\renewcommand{\labelenumi}{(\roman{enumi})}
In recording a spectrum, two separate independent steps are involved:
\begin{enumerate}
 \item sample preparation, and
 \item measurement of the spectral response of the sample.
\end{enumerate}
Here, equation~(\ref{eq:Bfeqg}) is a model for (ii). For (i), a 
separate model is required which, if convolutional,  would take the form
$$
f = S*\hat{f}
$$
where $\hat{f}$ denotes the ideal form of the sample and $S$ the smoothing
resulting from the sample preparation.

In monodisperse molecular weight distribution modelling in rheology,
$\hat{f}$ is a Dirac delta function of a single molecular weight~\cite{ThimmFMH2000}.
In NIR spectroscopy, $\hat{f}$ is a collection of particles of a single 
size and homogenously arranged~\cite{BurgerG2007}. In mass spectroscopic
peptide finger printing, $\hat{f}$ is a sum of Dirac delta function molecular
weights corresponding to a perfect segmentation of the protein into individual
molecular weight peptides.

For the analysis in the sequel, this relationship between $f$ and $\hat{f}$
 justifies the conclusion below that $f$
is an $L_2$ function. If it was assumed that $f$ was a Dirac delta function,
the analysis below would have to be performed with weaker norms.

The paper has been organized in the following manner.
The concept of a DHS is introduced in Section 2 as a special case of 
a variable Hilbert scale (VHS).
The generalization of B\'egout and Soria~\cite{BegS06} of the H\"older inequality 
for $L^p$ spaces is first used in Section 2 to derive a general form for 
variable Hilbert scales interpolatory inequalities,
which is then utilized to derive the basic DHS interpolation inequality. 
The application of DHS in the construction for error estimates for the approximate solution
of inverse operator equation problem is examined in Section 3
along with a discussion of analytic differentiation,
source conditions and the application of the results to deconvolutional
sharpening.

\section{Variable and dilational Hilbert scales}

The \emph{spectral theorem} for a positive definite selfadjoint operator $T$ on a Hilbert space $H$ 
(see, e.g.,~\cite{Rud73}) introduces (for each $T$) a family  $E(\lambda)$ of orthogonal projections 
such that $E_{f,g}(\lambda) := (E(\lambda) f, g)$ defines a Stieltjes measure on  $(0,\infty)$ and
$$
  (Tf, g) = \int_0^\infty \lambda\, d E_{f,g}(\lambda)
$$
for all $g\in H$ and $f\in D(T)$ (the domain of $T$). One can see that $(E(\lambda) f, g)$ is
right semi-continuous. 
In applications, the operator $T$ is often a differential operator like the 
Laplacian.
Intuitively, this representation generalizes the concept of a singular
value decomposition of a matrix.

Following the definition in~\cite{Heg95}, let any  measurable function 
$\phi:(0,\infty)\rightarrow(0,\infty)$ be called an \emph{index function}. 
Then
$$
  (f,g)_\phi :=  \int_0^\infty \phi(\lambda) dE_{f,g}(\lambda),\quad \|f\|^2_\phi=(f,f)_\phi
$$
is a densely defined bilinear form on $H$ with scalar product $(\cdot,\cdot)_\phi$.
Let the closure of the domain of this
bilinear form be denoted $H_\phi$. Note that $H_\phi$ then becomes a Hilbert space with scalar
product $(\cdot,\cdot)_\phi$ and we call the set of all possible $H_\phi$ \emph{a variable 
Hilbert scale}.

In~\cite{BegS06}, B\'egout and Soria introduce a generalisation of 
the H\"older inequality for $L^p$ spaces.
Here, we adapt this theorem to derive a key property of the above variable Hilbert scales. 
In the following, it will be assumed that 
products of index functions are defined point wise; i.e.,
$$
  (\phi \theta) (\lambda) = \phi(\lambda) \theta(\lambda).
$$
\begin{theorem}[Generalised H\"older inequality]
  Let $\phi,\psi$ and $\theta$ be three index functions and $\Phi$ and $\Psi$ be concave functions
$(0,\infty)\rightarrow (0,\infty)$ such that a.e.
\begin{equation}
  \label{eq:Holcond}
  1 \leq \Phi(\phi(\lambda)) \, \Psi(\psi(\lambda)).
\end{equation}
Then 
\begin{equation}
  \label{eq:genHol}
  1 \leq \Phi(\|f\|^2_{\phi\theta}/\|f\|_\theta^2)\, \Psi(\|f\|_{\psi\theta}^2/\|f\|_\theta^2)
\end{equation}
for all $f \in H_{(\phi+\psi+1)\theta}$.
\end{theorem}
\begin{proof}
  Let the measure $\nu$ be defined by 
$d\nu(\lambda) = \|f\|_\theta^{-2} \theta(\lambda) dE_{f,f}(\lambda)$
  for $f\in H_{(\phi+\psi+1)\theta}$. 
By definition, 
because $\theta$ is positive and the integral of $d\nu$ equals 1, 
$\nu$ is a probability measure and, by 
  equation~(\ref{eq:Holcond}), one obtains
$$
  1 = \int_0^\infty d\nu(\lambda)\leq \int_0^\infty \Phi^{1/2}(\phi(\lambda))\, 
\Psi^{1/2}(\psi(\lambda)) \, d\nu(\lambda).
$$
After taking the square of the right-hand side, 
an application of the Cauchy-Schwartz inequality yields
$$
  1 \leq \int_0^\infty \Phi(\phi(\lambda)) d\nu(\lambda) \, 
\int_0^\infty \Psi(\psi(\lambda)) d\nu(\lambda).
$$
Since both $\Phi$ and $\Psi$ are concave, one can use the (inverse) Jensen inequality to obtain
$$
  1 \leq  \Phi\left(\int_0^\infty \phi(\lambda) d\nu(\lambda)\right)\, 
          \Psi\left(\int_0^\infty \psi(\lambda) d\nu(\lambda)\right).
$$
The required inequality~(\ref{eq:genHol}) is then obtained by replacing $d\nu(\lambda)$ by its definition.
\end{proof}
A related inequality is obtained by choosing 
$\phi(\lambda)=\lambda^{-m}$, $\psi(\lambda)=\lambda^m$
and $\theta(\lambda)=\Phi(\lambda)=\Psi(\lambda)=1$. 
In fact, one recovers a special case of
a generalised H\"older inequality which holds for Sobolev spaces (see~\cite[p. 50]{Ada75}); 
namely,
$$
  \|f\|^2 \leq \|f\|_{-m} \, \|f\|_m
$$ 
where $\|f\|_{-m} := \|f\|_\phi$ and $\|f\|_m := \|f\|_\psi$.

The interpolation inequality for variable Hilbert scales (see theorem 2.2 in~\cite{Heg95}) 
can now be obtained as a direct consequence of the above generalised H\"older inequality:
\begin{corollary}[Interpolation inequality~\cite{Heg95}]
\label{cor2}
  Let $\phi,\psi$ and $\theta$ be index functions. 
  \begin{itemize}
 \item If $\phi\circ \psi^{-1}$ is concave, then
\begin{equation}
     \|f\|_{\phi\theta}^2 \leq \|f\|_\theta^2\; 
\phi\circ\psi^{-1} \left(\|f\|_{\psi\theta}^2/\|f\|^2_\theta\right),
      \quad f\in H_{(\phi+\psi+1)\theta}, \quad f \neq 0.
\end{equation}
 \item If  $\phi$ and $\psi$ are strictly increasing and $\psi\circ \phi^{-1}$ is convex, 
then
\begin{equation}
  \phi^{-1}\left(\frac{\|f\|^2_{\phi\theta}}{\|f\|^2_\theta}\right) \leq 
  \psi^{-1}\left(\frac{\|f\|^2_{\psi\theta}}{\|f\|^2_\theta}\right), 
   \quad f\in H_{(\phi+\psi+1)\theta} \quad f \neq 0.
\end{equation}
\end{itemize}
\end{corollary}
\begin{proof}  
Let $\Phi(\lambda)=1/\lambda$ and $\Psi(\lambda)=\phi\circ\psi^{-1}(\lambda)$. It follows that
$$
  \Phi(\phi(\lambda))\, \Psi(\psi(\lambda) = 
\frac{1}{\phi(\lambda)} \phi(\psi^{-1}(\psi(\lambda)))=1.
$$
Furthermore, $\Phi$ is concave, and, as $\psi\circ \phi^{-1}$ 
is convex and monotonically increasing, $\Psi$ is concave. 
The following inequality then follows from Theorem 1
$$
  \frac{\|f\|_{\phi\theta}}{\|f\|_\theta^2} 
\leq \phi(\psi^{-1}(\|f\|_{\psi\theta}^2/\|f\|_\theta))
$$
which can be rearranged to give the first inequality. 
The second inequality follows from the monotonicity of $\phi^{-1}$.
\end{proof}

The standard Hilbert scales inequality is recovered by choosing $\theta(\lambda)=\lambda^a$,
$\phi(\lambda)=\lambda^{b-a}$ and $\psi(\lambda)=\lambda^{c-a}$ with $b>a$ and $c>a$.
If $\sigma=(b-a)/(c-a) < 1$, then $\phi\circ\psi^{-1}(\lambda) = \lambda^\sigma$ is concave,
and the standard interpolation inequality for Hilbert scales is obtained
$$
  \|f\|_b \leq \|f\|_a^{1-\sigma} \|f\|_c^{\sigma},
$$
with $\|f\|_b:=\|f\|_\phi$, $\|f\|_a:=\|f\|_\theta$ and $\|f\|_c:=\|f\|_\psi$.
By choosing $\theta(\lambda)=1$ in the previous corollary, 
one obtains the version of the interpolation inequality which will be used in the sequel.
\begin{corollary}
  Let $\phi$ and $\psi$ be index functions such that $\phi\circ\psi^{-1}$ is concave. Then 
$$
  \|f\|_\phi^2 \leq \|f\|^2 \phi(\psi^{-1}(\|f\|_\psi^2/\|f\|^2)), \quad f\in H_\psi.
$$
\end{corollary}
A common situation in applications arises when $\phi(t)=t$ and $\psi(t)=t^m$ for some integer $m>1$.
For $m>1$, $\phi\circ\psi^{-1}(\lambda)=\lambda^{1/m}$ is concave, which, thereby, yields
$$
  \|f\|_1 \leq \|f\|^{1-1/m} \|f\|_m^{1/m}.
$$

\subsection{\label{sec:2.1} Dilational Variable Hilbert Scales (DHS)}

A new special family of variable Hilbert scales can be generated from 
a monotonically increasing index function 
$a: [0,\infty) \rightarrow [1,\infty)$ with $a(0)=1$, 
when the scales are all of the form $\phi(\lambda) = a(s\lambda)$, 
with $s>0$, and the corresponding norms are defined to be
$$
  \|f\|_s^2 = \int_0^\infty a(s\lambda) dE_{ff}(\lambda).
$$
The Hilbert space with norm $\|\cdot\|_s$ will be denoted by $H_s$. As the index functions are
obtained from dilations of the original function $a(\lambda)$, we will refer to this family of variable Hilbert
scales as \emph{dilational Hilbert scales} (DHS). 
Because the generating index function $a$ is 
monotonically increasing, it follows that, for $s\leq t$ and $f\in H_t$, $\|f\|_s \leq \|f\|_t$, which implies
the existence of a continuous embedding $H_s \hookrightarrow H_t$. Since $a(0)=1$, it follows that
$\|f\|_0 = \|f\|$ and $H_0=H$.

Though ordinary Hilbert scales (OHS) are not directly derivable as 
a special case from DHS, 
they do form a special subset of DHS as identified in
\begin{proposition}
  The (ordinary) Hilbert scales, defined by a self-adjoint operator $T\geq I$, are identical with the
   dilational variable Hilbert scales generated by $a(\lambda)=\exp(\lambda)$ and
 the operator $\log T$.
\end{proposition}
Here $I$ is the identity and $T \geq I$ means that the spectrum is contained in $[1,\infty)$. This
is required so that the logarithm $\log T$ defined by
$$
  \log T = \int_1^\infty \log(\lambda) dE(\lambda) = \int_0^\infty \mu\, dE(\exp(\mu))
$$
is well defined and positive definite. The proposition follows automatically on applying these 
special choices 
to $a(s\mu)=\exp(s\log \lambda) = \lambda^s$, where $\mu = \log(\lambda)$, and noting that the
spectral projections of $\log(T)$ are just $E(\exp(\mu))$.

The next theorem identifies an important subclass of DHS 
which genearte a quite special and useful set of interpolation inequality.
\begin{theorem}
  \label{thm:2}
   Let $\alpha(\lambda)$ satisfy the scaling relations
$$
  \alpha(\sigma \lambda) \leq \frac{1}{\sigma} \alpha(\lambda), \quad \sigma \in (0,1],
$$
and, for some $c\geq 0$, let the generating function of the DHS be
$$
  a(\lambda) = 1 + c\int_0^\lambda \exp\left(\int_1^t \alpha(s) ds\right) dt.
$$
Then one obtains, for every $\sigma \in [0,1]$, the interpolation inequality
$$
  \|f\|_{\sigma t}^2 \leq \|f\|^2 a(\sigma a^{-1}(\|f\|_t^2/\|f\|^2).
$$
\end{theorem}
\begin{proof}
   This theorem is an immediate consequence of corollary~\ref{cor2} with $\phi(\lambda)= a(\sigma \lambda)$
   and $\psi(\lambda) = a(\lambda)$, once it is established that the conditions of the corollary
   are fulfilled. 

   By definition, $a(\lambda)$ is differentiable with derivative
$$
  \dot{a}(\lambda) = c\exp\left(\int_1^\lambda \alpha(s) ds\right).
$$
   Because $c \geq 0$, it follows that $\dot{a}(\lambda) > 0$ for $\lambda > 0$ and, hence,
 that both
   $\phi$ and $\psi$ are strictly monotonically increasing.

   As $a$ is an integral with positive integrand, it follows that $\phi(\lambda) \leq \psi(\lambda)$.

   The remaining condition to show is the convexity of $\psi\circ \phi^{-1}$. Let
   $$
     \theta(\lambda) = \psi(\phi^{-1}(\lambda)) = a(a^{-1}(\lambda)/\sigma).
   $$
   From the chain rule, one obtains
$$
   \dot{\theta}(\lambda) = \frac{\dot{a}(a^{-1}(\lambda)/\sigma)}{\sigma \dot{a}(a^{-1}(\lambda))}.
$$
   The convexity of $\theta$ is equivalent with the monotonicity of $\dot{\theta}$, which, as $a$
   is monotonically increasing, is equivalent to the monotonicity of $\zeta(\lambda)=\dot{\theta}(a(\lambda))$.
   One now has
$$  
   \zeta(\lambda) = \frac{\dot{a}(\lambda/\sigma)}{\sigma \dot{a}(\lambda)}
$$
   and so, by the quotient rule of differentiation, 
$$
  \dot{\zeta}(\lambda) = \frac{\ddot{a}(\lambda/\sigma)}{\sigma^2\dot{a}(\lambda)} -
  \frac{\dot{a}(\lambda/\sigma) \ddot{a}(\lambda)}{\sigma \dot{a}(\lambda)^2}.
$$
Because $\dot{a}(\lambda) \geq 0$, the required monotonicity reduces to the scaling condition
$$
  \frac{\ddot{a}(\lambda/\sigma)}{\dot{a}(\lambda/\sigma)} \geq \sigma
  \frac{\ddot{a}(\lambda)}{\dot{a}(\lambda)}.
$$

From the definition of $a$, it follows that $\log \dot{a}(\lambda) = \log(c) + \int_1^\lambda \alpha(s) ds$ and,
hence,
$$
  \frac{\ddot{a}(\lambda)}{\dot{a}(\lambda)} = \frac{d}{d\lambda} \log \dot{a}(\lambda) = \alpha(\lambda).
$$
On replacing $\lambda$ with $\lambda/\sigma$ in the scaling condition for $\alpha$, the required
monotonicity and, hence, convexity are established.
\end{proof}

The simplest example arises on choosing $\alpha(\lambda)=1$ and $c=e$, which corresponds to taking
$a(\lambda)=\exp(\lambda)$. Another example, which is the limiting example for the convexity condition,
arises when $\alpha(\lambda)= (\gamma-1)/\lambda$ and $c=\gamma$,
which corresponds to taking $a(\lambda)= 1 + \lambda^\gamma$.
Even though this generates a one-parameter family of scales like 
the OHS family ,
it is a different family because its index function takes the form $a(s\lambda) = 1 + s^\gamma \lambda^\gamma$
where $s$ is the parameter indexing the Hilbert spaces $H_s$.



Remark. In the deliberations below, the interpolation inequalities derived above 
proved to be sufficient for the purposes required.
Though ordinary interpolation inequality can be used to bound the norms $\|f\|_{(1-\sigma)r+\sigma t}$
in terms of the norms $\|f\|_r$ and $\|f\|_t$, for any $0\leq r$, only the situation with
$r=0$ is discussed in the sequel.
Consequently, this represents a matter where the current theory might be generalised. 
In addition, other extensions might include
different choices for the functions $\Phi$ and $\Psi$ in the generalized H\"older inequality.

\section{\label{sec:3} Solving Ill-posed Problems}

Because of their special structure,
variable Hilbert scales interpolatory inequalities yield a natural
framework in which to derive upper bounds for the errors associated with 
approximate regularized solutions of ill-posed operator equations 
$Af=g$ for "non-exact data" $g_\epsilon$. 
For well-posed operator equations, they are not required. 
Assume, for simplicity, that $A$ is injective so that,
for every element $g$ in the range $R(A)$ of $A$, 
there is a unique $f$ such that $Af=g$. 
For well-posed $A$, $R(A)$ is the full Hilbert space $H$,
which guarantees that $g_\epsilon \in R(A)$ and,
furthermore, that the inverse $A^{-1}$ of $A$ exists and is bounded. 
Consequently, on assuming that $\|g_\epsilon - g\| \leq \epsilon$, 
one can immediately derive the upper bound
$$
 \|f_\epsilon - f\| \leq \|A^{-1}\|\, \|g_\epsilon - g \| \leq \|A^{-1}\| \,\epsilon
$$
without the need for an  interpolatory inequality.

For ill-posed problems, 
$R(A)$ is not closed and $A^{-1}$, if it exists, is not bounded. 
Consequently, unlike for well-posed problems,
there is no guarantee that there exists an $f_\epsilon$ such that the residual
$r = Af_\epsilon - g$ is equal to $g_\epsilon -g $.  
However, from the boundedness of $A$, it follows that 
$\|Af_\epsilon - g \| \leq \|A\| \|f_\epsilon -f\|$.
It therefore follows that small residuals is a \emph{necessary condition}
for the errors to be small. 
In a way, this condition corresponds to the \emph{consistency condition} 
in the Lax-Richtmyer equivalency theorem~\cite{Lin79}.
In fact, 
for the approximations $f_\epsilon$ generated by regularisation methods,
the pseudo-data $\tilde{g}_\epsilon := Af_\epsilon$
is a well-defined function of $f_\epsilon$. 
One can thereby interprete $f_\epsilon$ as the solution of 
$Af_\epsilon = \tilde{g_\epsilon}$. 
Consequently, 
the \emph{consistency error} corresponds to the difference 
$$
  g-\tilde{g}_\epsilon=A f - \tilde{g}_\epsilon = g - A f_\epsilon = -r,
$$
and, thereby, to the negative of the residual generated  by $f_\epsilon$.

However, because 
the available data are only known approximately as $g_\epsilon$, 
the value of the residual $Af_\epsilon-g$ cannot be determined and
utilized to assess the appropriateness of a given $f_\epsilon$.
A popular surrogate is the \emph{discrepancy} $Af_\epsilon - g_\epsilon$. 
The goal is to limit the choice of the $f_\epsilon$ so as to guarantee that 
$\|Af_\epsilon - g_\epsilon\| \leq \epsilon$. 
From the triangular inequality, it follows then that 
$$
   \|Af_\epsilon - g \| 
\leq \|Af_\epsilon - g_\epsilon\| + \|g_\epsilon - g\| \leq 2\epsilon,
$$
which is just twice the value obtained for a well-posed problem.

Interestingly, this bound for the residual,
associated with the reconstruction, 
is conceptually similar to the consistency bound used in the theory of the 
numerical solution of well-posed initial value problems and operator equations. 
In such situations, 
convergence is established by imposing a stability condition on the
numerical performance of the numerical method and then 
applying the Lax-Richtmyer equivalence theorem (see, e.g.,~\cite{Lin79}). 
For ill-posed problems, 
the interpolation inequality can be used to derive 
a similar (equivalence) theorem.
In fact, the \emph{source condition}, 
in the form  $f=S u$ with $\|u\| \leq C$ for some constant $C$, 
can be reinterpreted as a stability condition by assuming that
the reconstructions $f_\epsilon = S u_\epsilon$ are such that 
$\|u_\epsilon\| \leq C$ for some constant $C$.  

The operator $S$, which defines the source condition, 
is often chosen to take the form $S = (A^*A)^m$, $m>0$. 
On setting the generator $T$ of the Hilbert scales to be $(A^*A)^{-1}$,
the norms in which $f$ are bounded are thereby determined
along with the associated counterpart of equation (1). 
However, such source conditions,
which correspond to taking $\phi(t)=t^\alpha$ and $\psi(t)=t^\beta$,
do not always generate sharp enough error bounds.
This is resolved by defining $S$ to be some more general function of $A^*A$.
An example is given in Subsection 3.1,
where the differentiation of analytic functions is examined. 
It illustrates how variable Hilbert scales are able to yield substantially
improved error bounds. 

If it is assumed that there exists an operator $T$ and index functions
$\phi$ and $\psi$ such that $R(A) \subset H_\phi$ and $R(AT) \subset H_\psi$, 
one obtains a chain of three embedded Hilbert spaces 
$H_\psi \hookrightarrow H_\phi \hookrightarrow H$ with $g\in H_\psi$. 
Regularisation is then achieved by controlling the size
of $\|Af_\epsilon\|_\psi$. 

Here, it will be assumed that $\max (\|g\|_\psi, \|Af_\epsilon\|_\psi)
\leq C$ such that $\|r\|_\psi \leq 2C$.
One then obtains, for $e=f_\epsilon-f$, 
on using the corresponding interpolation inequality for variable Hilbert
scales (Hegland~\cite{Heg95}), 
the following bound
$$
  \|e\|^2 \leq \|r\|_\phi^2 \leq \|r\|^2 \phi\circ\psi^{-1}(\|r\|_\psi^2/\|r\|^2)
                       \leq \|r\|^2 \phi\circ\psi^{-1}(\|v\|^2/\|r\|^2)
$$
where $v=S^{-1} e=u_\epsilon - u$ satisfies the bound $\|v\|\leq \|u\| + C$.
 
If it is assumed that $\phi(t)/\psi(t)$ is monotonically decreasing such that 
$$
   \lim_{t\rightarrow 0} \phi(t)/\psi(t) = 0,
$$
a convergent rate can be derived.
A simple substitution shows that the right-hand side of the previous
bound is monotonically increasing in $\|r\|$. 
Consequently, if one invokes the discrepancy ansatz
(i.e., 
chooses $f_\epsilon$ such that $\|r\| \leq 2\epsilon$ and $\|r\|_\psi \leq 2C$),
it follows that
\begin{equation}
  \label{eq:errorbnd}
  \|e\| \leq 2\epsilon \sqrt{\phi\circ\psi^{-1}\left(C^2 \epsilon^{-2}\right)}
\end{equation}
which is decreasing monotonically to zero as $\epsilon \rightarrow 0$.

Consider, for example, the case where $\psi(t)=\exp(t)$ and $\phi(t)=t$.
Because $\psi(\phi^{-1}(t))=\exp(t)$ is convex and $t \leq \exp(t)$, 
the conditions given in Hegland~\cite{Heg95} are satisfied.
In this way, the following error bound holds
\begin{equation}
  \label{eq:errorbnd2}
  \|e\| \leq 2\sqrt{2(\log C/\epsilon)}\, \epsilon
\end{equation}
for sufficiently small $\epsilon$.

For DHS, the chain of three  Hilbert spaces becomes 
$H_t \hookrightarrow H_{\sigma t} \hookrightarrow H$. In addition, using
the interpolatory inequality for DHS of Theorem~\ref{thm:2}, one now
obtains, assuming $\|e\|\leq\|r\|_{\sigma t}$,
$$
  \|e\|^2\leq\|r\|_{\sigma t}^2 \leq \|r\|^2 a(\sigma a^{-1}(\|r\|^2_t/\|r\|^2)).
$$
For each of the examples considered at the end fo Section~\ref{sec:3.3}, it will
be shown that the assumed inequality $\|e\| \leq \|r\|_{\sigma t}$ is valid.

\subsection{Analytic differentiation exemplification}

Ordinary Hilbert scales correspond to choosing $\psi(\lambda)=\lambda^m$ 
and $\phi(\lambda)=\lambda$. Because $\phi\circ\psi^{-1}(\lambda) = \lambda^{1/m}$ is 
concave, it follows that
$$
  \|e\| \leq 2 C^{1/m} \epsilon^{1-1/m}
$$
which decays much slower to zero than the bound given by~(\ref{eq:errorbnd2}). 
Consequently, if it is applied to estimating the errors in the computation of
the derivative of observed data of an analytic function, suboptimal error bounds will result.

As exemplification, consider the determination of the derivative of a periodic function
which is the restriction of an analytic function onto the unit circle. 
In practice, this problem can arise when one 
needs to determine the electric field from measurements of an electrostatic potential. We
will consider a slightly simplified situation below, but the same arguments can
be applied to other applications and domains. Consider the unit complex sphere $B_1 \subset \C$.
Let $H^2(1)$ be the Hardy space of all analytic functions on $B_1$ with the scalar product
$$
  (f,g) = \frac{1}{2\pi} \int_0^{2\pi}\, \overline{f(e^{2 \pi i\theta)})}\, g(e^{2\pi i\theta}) d\theta,
  \quad f,g \in H^2(1), 
$$
and, for simplicity, we will also assume that $\int_0^{2\pi} f(e^{2\pi i \theta}) d\theta = 0$, the
general case can easily be derived from this case. 
There is a one-to-one correspondence between the potentials which are the real parts of the analytic
functions $f$ and $g$. It follows directly that the monomials form an orthonormal system of $H^2(1)$
and furthermore, that the norm in $H^2(1)$ is
$$
  \|f\|^2 = \sum_{k=1}^\infty |a_k|^2
$$
if $f(z) = \sum_{k=1}^\infty a_k z^k$ (recall that we assume that $a_0=0$).  Consider now the 
self-adjoint positive operator $T$ defined by
$$
  (Tf, g) = \sum_{k=1}^\infty k\, \overline{a_k} \, b_k.
$$
A differentiable function is one where the derivative is in $H^2(1)$. This requires that
$$
  \left\|\frac{dg(z)}{dz} \right\|^2 = \sum_{k=1}^\infty k^2 |b_k|^2 < \infty.
$$
If we set $H=H^2(1)$, then the differentiable functions are just the elements of the space $H_\phi$
defined previously using the operator $T$ above along with $\phi(t)=t^2$. 

Here, the source condition can be formulated as a condition on the (unperturbed) data $g$. 
If it is assumed that there are no singularities of the field in a sphere of radius $R > 1$, 
then the data $g(z)$ is analytic
in the sphere $B_R$ and in the Hardy space $H^2(R)$ of analytic functions on $B_R$ with a finite 
squared Hardy norm:
$$
  \|g\|^2 = \frac{1}{2\pi}\int_0^{2\pi} |g(Re^{2\pi i \theta})|^2\, d\theta
          = \sum_{k=1}^\infty R^{2k} |b_k|^2 < \infty, 
$$
where $g(z) = \sum_{k=1}^\infty b_k z^k$. It follows that $g\in H_\psi$ with $\psi(t)=R^{2t}$.
(Actually, it is the restriction of $g$ onto the unit sphere $B_1$ that is examined here.). 
Applying the VHS error bound~(\ref{eq:errorbnd}),  
one obtains the following bound for the differentiation of analytic functions
$$
  \|e\| \leq \frac{2\epsilon |\log(C/\epsilon)|}{\log R}.
$$
This yields an explicit demonstration that substantially better asymptotic
error bounds can be derived using VHS than can obtained using OHS. 

\subsection{Source conditions}

In some cases, the ``exact source condition'' that gives the maximum convergencen
rate, i.e., the maximal $\psi$ for which
$Af\in H_\psi$ may not be known. One approach, the Morozov method~\cite{Mor93} then stabilises
the norm $\|f_\epsilon\|$ of the solution. A more general method would lead to a bound on $\|r\|_\theta$
for some index function $\theta$. In this way, a chain of four Hilbert spaces
$$
   H_\psi \hookrightarrow H_\theta \hookrightarrow H_\phi \hookrightarrow H
$$
is generated.
The $\phi$ relates to the error norm which one would like to be small, the $\theta$ corresponds
to the norm which is controlled by the algorithm and the $\psi$ relates to the smoothness of the
solution. 
In the case discussed in the previous subsection for the embedding discussed in
the first part of Section~\ref{sec:3}, 
one has $\theta=\psi$ and, for the Morozov method, one has $\theta=\phi$.
For $\theta=\psi$, an 
algorithm based on $\|r\|_\theta$ would then yield the bound
$$
  \|e\|^2 \leq \|r\|^2 \phi \circ \theta^{-1} (\|r\|_\theta^2/\|r\|^2).
$$

Consider now the generalised Morozov method where $f_\epsilon$ is chosen so that 
$\|f_\epsilon\|_\theta \leq \|f\|_\theta$. 
(The Morozov method gets this bound by choosing $f_\epsilon$ to minimize the
$\theta$-norm.) 
Then one obtains 
\begin{align*}
   \|r\|_\theta^2  = & \|Af\|_\theta^2 - 2(Af, Af_\epsilon)_\theta + \|A f_\epsilon \|_\theta^2 \\
                   \leq & 2(g, r)_\theta \\
                   \leq & 2\|g\|_\psi \, \|r\|_{\theta^2/\psi}, 
\end{align*}
where the following variant of the Cauchy-Schwartz inequality~\cite{Heg95}
 has been used:
$$
  (g,r)_\theta \leq \|g\|_\psi\, \|r\|_{\theta^2/\psi}.
$$
In the classical Morozov situation, where $\theta=\phi$ and $\psi=\phi^2$, 
the following well-known $\sqrt{\epsilon}$ convergence estimate holds:
$$
  \|e\| \leq 2 \sqrt{\|g\|_\psi \epsilon}.
$$
Combining this with the corresponding interpolation inequality, 
error rates are derived for situations where
$\phi < \theta < \psi$.

\subsection{Application - Deconconvolution by sharpening\label{sec:3.3}}

The theory discussed above provides a natural framework in which 
to analyse the sharpening of broadened and possibly overlapping
spectroscopic peaks by deconvolution.  Because, as shown in the 
introduction, the convolution of two Gaussian peaks is a Gaussian 
peak, an OHS analysis can be applied successfully when the peak,
to be measured, and its broadening, that occurs as the result of its 
measurement, are both modelled as Gaussians.
In the DHS framework, general dilational-parameterized peaks can be 
analysed with similar facility. The VHS framework could be applied, 
but additional assumptions would have to be invoked,
like the ones given above in Theorem 2 for DHS, before error estimates 
as sharp and useful as those given below could be derived.
In essence, Theorem 2 generates a framework which allows the VHS
methodology, as a DHS methodology,
 to be applied directly to deconvolution by sharpening.

Spectroscopy reveals information about the chemical composition of samples 
and is an important tool in chemistry, physics, biology, astronomy
and related industrial applications. The data consist of a  
superposition of ``peaks''. In the case of overlapping peaks, 
their separation and identification poses a substantial challenge. 
Methods for performing such tasks are discussed in~\cite{HegA05}.
They have wide applicability and can also be used for deblurring
in image processing.

The widening of the peaks in a spectrum results from a ``diffusion" of 
information into neighboring frequencies. If this ``diffusion" is 
independent of location, it can be modelled as a convolution. 
In an $L_2(\R)$ Hilbert space context, the theoretical model takes the 
form
$$
 g(x) = B_\gamma*f(x) = \int_{-\infty}^\infty B_\gamma(x-y) f(y) dy, 
\quad B_\gamma \in L_2(\R),
$$
where $f\in L_2(\R)$ is the actual spectrum being measured, $g\in L_2(\R)$ 
is its measurement and $B_\gamma(x)=\gamma^{-1}B(x/\gamma)$ models the broadening 
that the measurement process has performed.
The observed spectral data $g_\epsilon$ are perturbed by ``observational noise''
such that $\|g_\epsilon - g \| \leq \epsilon$. The problem of recovering $f$ 
from $g_\epsilon$ is called \emph{deconvolution}. It can be severely ill-posed 
for smooth $B$. Because of the importance of spectroscopy across a wide range 
of applications, deconvolution is a prominent example of an ill-posed problem.

Rather than attempting to accurately perform the full deconvolution, it is more
sensible to
 ``sharpen'' the spectrum so that a  better identifiability of
the locations and number of peaks is achieved compared with that available 
from a visual inspection of the available data $g_\epsilon$. As explained
 in the Introduction,
this corresponds to finding the solution $z$ of 
$$
  B_{\gamma,\beta}*z = g
$$
for the data $g_\epsilon$ where the regularity of $z$ is determined by the
 \emph{source condition} 
$$
  z = B_\beta*f, \quad f\in L_2(\R).
$$
When applying the DHS interpolatory inequality in order to obtain bounds on the
error, it is appropriate at this stage to recast the source condition as a
condition on the data
$$
  g = B_\gamma*f, \quad, f\in L_2(\R).
$$
In order to utilise the DHS interpolation inequality,
it is necessary to introduce a DHS (as introduced in section~\ref{sec:2.1}) 
 which guarantees that $g$ is in $H_s$. One
way to achieve this is to choose a DHS such that
\begin{equation}
 \label{eq:source}
\|B_\gamma*f\|_s = \|f\|, \quad f\in L_2(\R)
\end{equation}
holds for some fixed $s$.
We now proceed to derive a DHS for which this is satisfied.

For convolutions, it is natural to generate the DHS
with $T=-d^2/dx^2$. Using the Fourier transform $\hat{f}$, one obtains
$$
 \int_0^\infty \lambda dE_{ff}(\lambda) = \frac{1}{2\pi}\int_\R
 \omega^2|\hat{f}(\omega)|^2\,d\omega = \frac{1}{2\pi}\int_0^\infty
 \omega^2 \left(|\hat{f}(\omega)|^2+|\hat{f}(-\omega)|^2\right)\, d\omega
$$
with $\lambda=\omega^2$ and 
$dE_{ff}(\lambda) = \frac{1}{2\pi}\left(|\hat{f}(\omega)|^2+|\hat{f}(-\omega)|^2\right)$.
The resulting DHS norm is given by
$$
  \|f\|_s^2 = \frac{1}{2\pi}\int_0^\infty a(s\omega^2)
\left(|\hat{f}(\omega)|^2+|\hat{f}(-\omega)|^2\right)\, d\omega.
$$

For the sequel, assumptions must be invoked regarding the Fourier
transforms $\hat{B}(\omega)$ of the peaks. First, it is assumed that the 
absolute values of the Fourier transforms are symmetric; i.e., that
$$
  |\hat{B}(\omega)| = |\hat{B}(-\omega)|.
$$
A straight forward application of the definition of $B_\gamma(t)$
and the Fourier transform
proves that $\hat{B}_\gamma(\omega) = \hat{B}(\gamma\omega)$. With this
and the assumed symmetry, one then obtains for the left hand side of 
condition~(\ref{eq:source})
$$
 \|B_\gamma*f\|_s^2 = \frac{1}{2\pi}\int_0^\infty a(s\omega^2) 
  |\hat{B}(\gamma\omega)|^2
  \left(|\hat{f}(\omega)|^2+|\hat{f}(-\omega)|^2\right)\, d\omega.
$$
Since
$$
  \|f\|^2 = \frac{1}{2\pi}\int_0^\infty
  \left(|\hat{f}(\omega)|^2+|\hat{f}(-\omega)|^2\right)\, d\omega,
$$
the required condition~(\ref{eq:source}) holds if $s = \gamma^2$ and 
\begin{equation}\
 \label{eq:a}
  a(\lambda) = \frac{1}{|\hat{B}(\sqrt{\lambda})|^2}.
\end{equation}

We now show that, for  three of the most important examples of peaks, this
choice does indeed lead to a DHS for which Theorem~\ref{thm:2} can be
applied. For Gaussian peaks, one has $\hat{B}(\omega)=e^{-\omega^2/2}$
and, by equation~(\ref{eq:a}), it follows for this case that 
$a(\lambda) = \exp(\lambda)$. For this case, the DHS coincides
with the OHS. For the second example, consider exponential peaks for which
 $\hat{B}(\omega) = 1/(1+\omega^2)$. In this situation, 
$a(\lambda) = (1+\lambda)^2$. Going back to the definition of 
the function $\alpha(\lambda)$ which generates $a$, one finds that
$\alpha(\lambda) = (1+\lambda)^{-1}$. It follows that
$\alpha(\sigma \lambda) \leq \alpha(\lambda)/\sigma$. Consequently, such
peaks also generate 
a DHS for which the interpolation inequality holds. 
For the third class, $\hat{B}(\omega) = e^{-|\omega|}$.
It follows that $a(\lambda) = \exp(2\sqrt{\lambda})$ and hence 
$\alpha(\lambda) = (1-1/(2\sqrt{\lambda}))/\sqrt{\lambda}$ for which
the inequality $\alpha(\sigma\lambda) \leq \alpha(\lambda)/\sigma$ 
also holds. Thus, the  corresponding DHS interpolation theorem
holds for the members of this class.

We have thus seen that, for three of the most important classes of spectra, 
the interpolation theorem holds for the DHS generated. However, the
theorem only provides a bound for the norm $\|r\|_{s\sigma}$ of the residual
$r$. 
It is now necessary to use it to 
 bound  the error norm $\|e\|$.
Initially, a condition is imposed on the Fourier transform of the peaks  
which guarantees the bound
$$
  \|e\| \leq \|r\|_{s\sigma}.
$$
Since $r=B_{\gamma,\beta}*e$, it can be shown that, under
appropriate conditions, 
$$
  \|z\| \leq \|B_{\gamma,\beta}*z\|_{s\sigma}, \quad z\in L_2(\R).
$$
As the norm is a continuous function, it suffices to show that this
holds for a dense subset in $L_2(\R)$.
Furthermore, since $z=B_\beta*f$,
it is only necessary to prove that 
$$
  \|B_\beta*f\| \leq \|B_\gamma*f\|_{s\sigma}^2.
$$
From the Fourier transform, it is clear that this inequality
holds when 
$$
  |\hat{B}(\beta\omega)|^2 \leq \frac{|\hat{B}(\gamma\omega)|^2}%
                         {|\hat{B}(\gamma\sqrt{\sigma}\omega)|^2},
  \quad \omega \in \R.
$$
Some rescalings reveal that there exists a $\sigma < 1$ for which
this holds if and only if, for each scalar $0<\tau<1$, there exists
a $0<\rho<1$ such that
$$
  |\hat{B}(\tau \omega) \hat{B}(\rho\omega)| \leq 
  |\hat{B}(\omega)|, \quad \omega \in \R,
$$
or, in terms of the functions $a$ defining the DHS, the condition
becomes: for every $\tau$, there has to exist a $\rho$ such that 
\begin{equation}
  \label{eq:cond9}
  a(\lambda) \leq a(\tau\lambda)\, a(\rho \lambda), \quad \lambda > 0.
\end{equation}
In the first example above, since $a(\lambda)=e^\lambda$, 
equality in~(\ref{eq:cond9}) 
is guaranteed if $\tau+\rho=1$. In the second example, since 
$a(\lambda)=(1+\lambda)^2$ and, because $(1+\tau\omega)(1+\rho\omega)
\geq 1+ (\tau+\rho)\omega$, the inequality holds whenever
$\tau+\rho=1$. In the last example, since $a(\lambda) = 
\exp(2\sqrt{\lambda})$, equality holds if  
 $\sqrt{\tau}+\sqrt{\rho}=1$. It follows that, 
for the  peaks  considered, the choice of
DHS suggested above does allow for the application of the interpolation
inequality. In addition,  the interpolation inequality 
provides a bound for the error associated with the reconstruction procedure.
Of course, the peak sharpening procedure does
need to satisfy the conditions discussed in the previous sections.

For the Gaussian, exponential and rational peaks, the
conditions required for the application of the DHS interpolation inequality
are satisfied. For the resulting bounds on the norm of $e=z_\epsilon-z$,  
one obtains
\begin{enumerate}
 \item the usual OHS error bound for the case of Gaussian peaks
       where $a(\lambda)=e^\lambda$:
$$
  \|e\| \leq \|r\|_s^\sigma \|r\|^{1-\sigma};
$$
 \item a convex combination for the case of exponential peaks
   where $a(\lambda)=(1+\lambda)^2$:
$$
  \|e\| \leq \sigma \|r\|_s + (1-\sigma)\|r\|;
$$
and 
 \item for the rational peak where $a(\lambda)=\exp(2\sqrt{\lambda})$:
$$
 \|e\| \leq \|r\|_s^{\sqrt{\sigma}} \|r\|^{1-\sqrt{\sigma}}.
$$
\end{enumerate}
In summary, one obtains from the interpolation inequality
three types of ``convex combinations'' for
the errors, with the  $\sigma$ determined 
by the reconstruction method. It follows from what was said above that
$\sigma$ satisfies $\sigma + \beta/\gamma=1$ for the first two shapes 
and $\sqrt{\sigma} + \sqrt{\beta/\gamma}=1$ for the third. 
By tuning $\beta$, one thereby tunes the parameter $\sigma$.  
Ideally, one would like to choose $\sigma=0$ to get the smallest
possible error. This, however, corresponds to the case $\beta=\gamma$;
i.e., to the case where no sharpening is done. Conversely, $\beta=0$
corresponds to ``full sharpening'' but in such situations, the errors will be 
large. In the practical application of the above bounds, it is necessary 
to carefully consider how the trade-off is performed between the
amount of sharpening and error in
the reconstruction of the sharpened peaks.

\section*{Acknowledgement}
Both authors wish to acknowledge the financial support received from the
Radon Institute of Computational and Applied Mathematics to participate in
the Special Semester on 
"Quantitative Biology analyzed by Mathematical Methods". 
Among other things, it gave them the opportunity to finalize this
paper during their visit.

\bibliographystyle{plain}
\bibliography{paper}
 
\end{document}